\documentclass[12pt]{article}
\usepackage{amsmath,amsfonts,amssymb,latexsym}
\usepackage[all]{xy}
\usepackage{theorem}
\newcommand{\CC}{{\mathbb C}}

\newcommand{\ZZ}{{\mathbb Z}}

\def\O{{\cal O}}
\def\F{{\cal F}}
\def\C{{\cal C}}
\def\pp{\:\raisebox{-0.1ex}{$\stackrel{p}{\to}$}\:}
\def\qq{\:\raisebox{-0.1ex}{$\stackrel{q}{\to}$}\:}

\newcommand{\G}{{\frak g}}

\newtheorem{Theorem}{Theorem}[section]
\newtheorem{proposition}[Theorem]{Proposition}
\newtheorem{Lemma}[Theorem]{Lemma}

\theoremstyle{definition}

\theorembodyfont{\upshape}
\newtheorem{Remark}[Theorem]{Remark}
\newtheorem{Proof}{Proof}

\newtheorem{remark}{Remark}

\newtheorem{definition}{Definition}

\title{Torsors on Affine Varieties}
\author{S. Subramanian}
\date{}
\begin{document}
\maketitle

\begin{abstract}
Let $X$ be a smooth affine algebraic variety over the field of complex numbers which is contractible. Then every algebraic $G$-torsor on $X$ is algebraically trivial if $G$ is a semi-simple algebraic group. We also show that if $X$ is a smooth affine algebraic variety such that $\Omega^1_X$ is trivial and $X$ is topologically simply connected with $H^i(X,\ZZ)=0$ for $i=1,2$ and $3$, then every algebraic $G$-torsor on $X$ is algebraically trivial for a semi-simple algebraic group $G$.  
\end{abstract}

\section{Introduction} 

Let $k$ be an algebraically closed field of 
characteristic $0$ contained in  $\CC$. Let $X$ be a smooth affine 
variety over $k$ and let $G$ be an affine algebraic 
semisimple  group over $k$. We use the words torsor and principal 
$G$-bundle  interchangeably. The question we are interested in is, when 
are $G$  bundles on $X$ trivial? We let $X_\CC$ denote the base change 
of 
$X$ to $\CC$, the field of complex numbers.  In this article, we show the 
following

\paragraph*{Main Theorem} Let the singular cohomology  of $X_\CC$ with 
integer coefficients be trivial, let $\Omega^1_X$ be algebraically
trivial, and assume that the topological 
fundamental group of $X_\CC$ is trivial.  Then every algebraic $G$-bundle
on $X$ is algebraically trivial.
(See Theorem 6.1 , Remark 6.2 and Remark 6.4)

\begin{Remark}
For the definition and basic properties of torsors,
see \cite{M}, Ch.3, section 4.
In our case, $G$-torsors are locally isotrivial, i.e., every point
has a finite etale neighborhood over which the torsor is trivial (this is
a theorem of Grothendieck, see \cite{Ra}, Lemme XIV 1.4).
The main theorem states that under the hypotheses mentioned in the theorem,
every algebraic $G$-bundle is globally algebraically trivial in the
Zariski topology.
\end{Remark}

\begin{Remark}
The theorem of Quillen-Suslin (Serre's conjecture)
states that every projective module over a polynomial ring over $\CC$
is free. The case where the module has trivial determinant is a case of the
above theorem. It is not difficult to see that rank one projective
modules are free on a polynomial ring over $\CC$. See \cite{L}.
\end{Remark}

\section{A Reduction}

  In this short section, we will show that in order to prove the main
theorem over an algebraically closed field $k$ contained in the field of
complex numbers, it is enough to consider the field of complex numbers.
We have the following

\begin{Lemma}Let $k$ be an algebraically closed field contained in
$\CC$ and let $X$ be an affine variety defined over $k$. Let $H$ be an
affine algebraic group of finite type defined over $k$. Let $E_H \to X$
be a $H$-torsor on $X$ such that $E_H \otimes \CC$ is trivial over
$X_{\CC}$. Then $E_H$ is trivial on $X$ over $k$.
\end{Lemma}

\begin{Proof} Since $E_H \otimes \CC$ is trivial over
$X_{\CC}$, we can assume that there is a finitely generated $k$-algebra
$B$ such that $E_H \otimes B$ is trivial over $X_B$,where $X_B$ denotes
the base change of $X$ to $B$. Since $k$ is algebraically
closed we can find a closed  point of Spec $B$ which
is rational over $k$, i.e., surjection
$$B \to k$$
such that the composite
$$k \to B \to k$$
is the identity. We get
$$X \pp X \otimes_k B \qq X $$
with the composite identity. Since $E_H\otimes B$ is trivial over $X_B$,
it follows that $q^*E_H$ is trivial over $X_B$. Hence $p^*q^*E_H$ is trivial
over $X$. Therefore $E_H$ is trivial on $X$ over $k$.  \hfill{Q.E.D.}
\end{Proof}

\section{The Differential Complex}

Let $k$ be an algebraically closed field of characteristic $0$.  Let 
$X$ be a smooth, affine  variety over $k$.  Let $G$ be 
an affine algebraic group over $k$ of finite type.  Let $\pi: E 
\rightarrow X$ be a principal $G$-bundle.

Let $\Omega^i_E$ denote the bundle of $i$-forms on $E$ and let 
$\Omega^i_X$ denote the bundle of $i$-forms on $X$.  We have the natural 
morphism.
$$
\pi^*\Omega^i_X \rightarrow \Omega^i_E \ \ \ \ i=1,\ldots, n.
$$
Since the group $G$ acts on $E$ (on the right, by convention), $G$ also 
acts on $\Omega^i_E$. We can consider  the sheaf of germs of invariant 
$i$-forms for the $G$-action.  This sheaf descends to a vector bundle 
$A^i(E) $ on $X$.

Further, the canonical map
$$\pi^*\Omega^i_X \rightarrow \Omega^i_E$$
induces an inclusion of  bundles
$$
0\rightarrow \Omega^i_X\rightarrow A^{i}(E).
$$
We now consider the exterior differential 
$$
d_E : \Omega^i_E \rightarrow \Omega^{i+1}_E, \ \ i=0,1,\ldots
$$
on $E$.  We observe that $d_E$ is equivariant for the $G$-action and 
hence descends to an operator (also denoted by $d_E$)
$$
d_E: A^i (E) \rightarrow A^{i+1} (E) \ \ \ i=0,1,\ldots
$$
We further obtain a commutative diagram 
$$
\xymatrix{
A^i(E) \ar[rr]^{d_E}&&   A^{i+1}(E)   \\
 \ar[u] \Omega_X^i  \ar[rr]^{d_X}&&\Omega_X^{i+1}  \ar[u] 
}
$$
where $d_X$ denotes the exterior differential on $X$, for 
$i=0,1,\ldots,\ (\Omega_X^{i}$ is taken to be $\O_X$, the 
structure sheaf, for $i=0$).

Since $X$ is affine, there is a splitting of the inclusion of bundles 
$$
\Omega^i_X \subset A^i(E)
$$
for $i=1,\ldots, $ on $X$. This follows from the fact that on an
affine variety,a short exact sequence of vector bundles splits as
$H^1(X,F)$ vanishes for any vector bundle $F$.

We now recall Grothendieck's algebraic de Rham theorem(see \cite{GH},
page 453);
it says that for 
any smooth affine algebraic variety over the field $\CC$ of complex 
numbers, the cohomology of the algebraic  de Rham complex of the variety 
is isomorphic to the singular cohomology of the variety with complex 
coefficients.  In particular, if we assume that the singular cohomology 
with complex coefficients of the variety is trival, it implies that the 
algebraic de Rham complex is globally exact.\\

We now have 

\begin{proposition}Let $X$ be a smooth affine variety over 
$\CC$ such that $X_\CC$ has trival singular cohomology with complex 
coefficients, and let $\Omega^1_X$ be algebraically trivial.
Let $G$ be semisimple and simply connected. Then these is 
a system of splittings
$$
s_i : A^i(E) \rightarrow \Omega^i_X, i = 1, 2,  
$$
such that
$$
d_X \circ s_1 = s_2 \circ d_E.
$$
\end{proposition}

\begin{Proof}
By Grothendieck's algebraic de Rham theorem (see 
paragraph above), the algebraic de Rham complex of $X_\CC$ is exact 
globally.  However, we observe that if a $i$-form on $X$ defined over 
$k$ is exact over $\CC$, then it is exact over $k$.  It follows that the 
algebraic de Rham complex of $X$ over $k$ is globally exact.

We first note that $A^0(E)=\O_X$ and $\Omega^0_X=\O_X$.  So 
to find the system of splittings, we start with $A^1(E)$.  Since $X$ is 
affine, there exists a splitting $s_1: A^1(E)\rightarrow \Omega^1_X$ of 
the canonical inclusion,  $\Omega^1_X \subset A^1(E)$
This gives a direct sum decomposition
$$
A^1(E) =\Omega^1_X \oplus B^1_E
$$
where $B^1_E$ is isomorphic to the quotient of $A^1(E)$ by $\Omega^1_X$.

To find $s_2$, we proceed as follows:

Let the splitting $s_1$ define the direct sum decomposition
$$
A^1 (E) =\Omega^1_X \oplus B^1_E.
$$
We now consider $d_E(B^1_E) \ \subset\  A^2(E)$.

We now claim that 
$$
d_E(B^1_E)\  \cap \ \Omega^2_X=0.
$$
For otherwise, let $\omega$ be a $2$ form on $X$ in $d_E(B^1_E)$.  We 
know that $d_X\omega = d_E \omega$ if $\omega$ is considered as an 
element of $A^2(E)$.  Since $\omega \in \ d_E(B^1_E)$, and $d_E^2=0$, it 
follows that
 $d_X \ \omega=0$.  By the exactness of the algebraic de Rham complex, 
it follows that there is a $1$-form $\eta$ on $X$ such that
$\omega =d_X \eta$(recall that $\omega \in \ d_E(B^1_E)$).  
It follows that there exists $\eta' \ \in \ B^1_E$ such that $d_E \ 
\eta' = d_E \ \eta$.
\end{Proof}

We need the following two lemmas. We will resume the proof of the 
proposition after the lemmas:

\begin{Lemma}
Let $p: Y_1 \rightarrow Y_2$ be a finite Galois etale
morphism of schemes over a field of characteristic zero such that 
$Y_1$ is affine.  Then $Y_2$ is affine.
\end{Lemma}

\begin{Proof} See Exercise 3.19, page 94, in Hartshorne's book \cite{H}. 
\end{Proof}

\begin{remark} It is a theorem of Chevalley that the image of a
finite surjective morphism from an affine variety is affine.
\end{remark}

\begin{Lemma} 
The total space $E$ is affine.
\end{Lemma}

\begin{Proof} We know that $\pi: E\rightarrow X$ is locally 
trivial in the etale topology on $X$.  Since the structure group $G$ is 
affine, this implies by the above lemma (Lemma 3.2), that every point 
$x$ in $X$ has a Zariski neighourbood $U$ such that $\pi^{-1}(U)$ is affine.

Let $\F$ be any  coherent sheaf on the total space $E$.  By 
Proposition 8.1,  page 250 in Hartshorne's Algebraic Geometry (see \cite{H}), it follows that the higher direct images $R^i\pi_* \F, \ i\geq 1$  are zero 
(by the remark in the above paragraph).  This implies that 
$H^i(E,\F)\cong \ H^i(X, \pi_*\F)$  for $i\geq 0$. However, $X$ is 
affine, and $\pi_*\F$ is quasi coherent on $X$, so 
$$
H^i(X, \pi_* \F)=0\ \mbox{for} \ i\geq 1.
$$
Hence $H^i(E,\F)=0$ for   $i\geq 1$ and hence $E$ is affine. 
\hfill{QED}.
\end{Proof}

\begin{remark}We observe that to show that $E$ is affine, we
can argue as follows:since $E \to X$ is a $G$-torsor, there is a 
free action of $G$ on $E$, whose quotient map is the given map
$E \to X$. By Proposition 0.7, page 14 in \cite{Mu}, the map $E \to X$
is affine. Since $X$ is affine, $E$ is affine.
\end{remark}

\begin{remark}
We observe that if $\pi: E\to X$ is a $G$-torsor, then usually the structure morphism $\pi$ is assumed to be an affine morphism. 
\end{remark}
 
 We now resume our proof of Proposition 3.1.  If $G$ is simply 
connected, and 
 under our hypothesis that $H^1(X, \ZZ)=0$, it follows that the first 
 singular cohomology of $E$ is zero. This can be seen as follows: the
homotopy exact sequence and the hypothesis that $G$ is simply
connected implies that the topological fundamental group of $E$
is isomorphic to that of $X$; therefore their abelianisations are
isomorphic and so the vanishing of $H^1(X, \ZZ)$ implies the
vanishing of $H^1(E, \ZZ)$ (a proof of the justification of using
the homotopy exact sequence for an algebraic torsor is given in
the proof of Theorem 5.3 below).
 Since $E$ is affine (Lemma 3.3 above), by Grothendieck's algebraic de Rham theorem it follows that 
 if $\psi$ is a 1-form on $E$ such that $d_E \ \psi=0$, then there is a 
 function $\varphi$ on $E$ such that $d_E \ \varphi =\psi$.  Now $d_E \ 
 \eta'=d_E\ \eta \Rightarrow d_E (\eta-\eta')=0$, and hence there is a 
 function $\varphi$ on $E$ such that $d_E \ \varphi =\eta-\eta'$.  
However,  $\eta-\eta'$ is $G$-invariant, and hence $d_E \varphi$ is 
$G$-invariant.  For $g \ \in \ G, g^*\ d_E \ \varphi =d_E \ g^* \ 
\varphi = d_E \ \varphi$ and hence 
$$d_E (g^* \ \varphi-\varphi)=0
$$
$\Rightarrow \ g^* \ \varphi - \ \varphi=$ constant $C_g \ \in \ k$,
 since $E$ is smooth and connected (as $X$ and $G$ are smooth and connected), and so the zeroth de Rham cohomology of $E$ is $k$.  

We observe that this defines a homomorphism from $G$ to the additive 
group (which is abelian).  Since we assume that $G$ is semisimple, this 
homomorphism is trivial, and hence we can choose a $G$-invariant 
function $\varphi$ such that,
$$
d_E \ \varphi =\eta \ - \ \eta'
$$
Hence 
$$
\eta' = \eta - \ d_E \ \varphi
$$
From the exact sequence 
$$
0\rightarrow \Omega^1_X \rightarrow  A^1(E) \rightarrow 
A^1(E)/\Omega^1_X\rightarrow 0
$$
and since both $\eta$ and $d_E \ \varphi$ are in $\Omega^1_X$, it 
follows 
that the image of $\eta - \ d_E \ \varphi$ in $A^1(E)/\Omega^1_X$ is 
zero. Since $\eta' \ \in \ B^1_E$, it follows that $\eta'=0$.  Hence
$\omega = d\eta'=0$, and we have proved the claim that 
$$
d_E (B^1_E) \cap \ \Omega^2_X=0
$$
We observe that $d_E (B^1_E)$ denotes the subbundle of $A^2(E)$
generated by the sections $d_E\phi$ of $A^2(E)$ as $\phi$ varies over
$H^0(X,B^1_E)$. This can be done as $d_EH^0(X,B^1_E)$ is a $k$-vector
subspace of $H^0(X,A^2(E))$ since $d_E$ is $k$-linear
(see Appendix 1, Proposition A.1: the appendix shows that the sections $d_E \phi$ generate a subbundle of $A^2(E)$ and not just a subsheaf).
We now have, on $X$
$$
0\rightarrow \Omega^2_X \oplus d_E(B^1_E) \rightarrow A^2(E) \rightarrow Q^2
\rightarrow 0
$$
We choose a splitting of $$A^2(E)\rightarrow Q^2 \rightarrow 0$$(which
exists since every short exact sequence of vector bundles on $X$ is
split, as observed before)
to get a 
direct sum decomposition
$$
A^2(E)=\Omega^2_X\oplus d_E(B^1_E) \oplus \stackrel{\sim}{Q^2}
$$
and this defines a splitting
$$
s_2 : A^2(E) \rightarrow \Omega_X^2
$$
satisfying
$$
d_X \ s_1 = s_2 d_E 
$$
\hfill{QED}

\section{Connections on induced vector bundles}

 Let $X$ be a smooth affine 
variety over $k$ such that $X_\CC$, the base change of $X$ to $\CC$, has
trivial singular cohomology. Let $\Omega^1_X$ be algebraically trivial.
Let $G$ be a semi-simple, simply connected group over $k$ and let 
$\pi : E\rightarrow X$ be a $G$ bundle on $X$.  Let $V$ be a rational $G$-module over 
$k$.  
Then we can form the associated vector bundle $\stackrel{\sim}{V} =E(V)$
on $X$.\\

On $E$, we have the complex 
$$
V \otimes \O_E \stackrel{d_E}{\rightarrow} (V\otimes \O_E) \otimes \Omega^1_E 
\stackrel{d_E}{\rightarrow} (V \otimes \O_E) \otimes \Omega^2_E
$$
where $V \otimes \O_E$ is the trivial vector bundle with fibre the vector space $V$ on which 
$G$ acts, and $d_E$ is the exterior differential operator on $E$.  Since,
as we have 
observed before, $d_E$ is equivariant for the $G$ action, we can go modulo the $G$ action, 
and we obtain,
$$
\stackrel{\sim}{V} \ \stackrel{d_{\stackrel{\sim}{V}}}{\rightarrow}\  
\stackrel{\sim}{V} \otimes A^1 (E)\stackrel{d_{\stackrel{\sim}{V}}}
{\rightarrow}\  \stackrel{\sim}{V} \otimes A^2(E)
$$
on $X$. Since $d^2_E=0$, we  also have $d_{\widetilde{V}}^2 =0$.

By Proposition 3.1, we have splittings 
$$
s_i : A^i(E) \to \Omega_X^i  \ \ \ \ \ i=1,2
$$
such that $d_X \ s_1 = s_2 \ d_E$.  We thus obtain
$$
\stackrel{\sim}{V}\stackrel{\nabla}{\rightarrow} \stackrel{\sim}{V}\otimes 
\Omega^1_X \stackrel{\nabla}{\rightarrow} \stackrel{\sim}{V}\otimes \Omega^2_X
$$
on $X$, where $\nabla$ is a connection on $\stackrel{\sim}{V}$ such that
$\nabla^2=0$.  In other words, we have constructed a flat connection on 
$\stackrel{\sim}{V}$.

\begin{Remark}We observe that $s_1$ (as in Proposition 3.1)
defines a connection on $E$. However, $\nabla$ is not this connection,
though $s_1$ is used to define $\nabla$.
\end{Remark}

\section{The Monodromy}

Let $X$ be a smooth affine variety over $k$.  We assume that the singular cohomology 
of $X_\CC$ with integer coefficients is trivial (i.e., $H^0(X_\CC, \ZZ) =\ZZ$, $H^i(X_\CC, \ZZ) =0$ for $ i\geq 1$). Let $\Omega^1_X$ be algebraically
trivial. 
Let $G$ be a semi simple, simply connected algebraic group over $k$.  Let $\pi:E
\rightarrow X$ be a $G$-bundle on $X$.  By Proposition 3.1 above, there are splittings 
$s_i : A^i(E)\rightarrow \Omega^i_X \ \ \  i=1,2$, such that $d_X \circ 
s_1 =s_2 \circ d_E$.\\

We fix one such set of splittings to begin with. Given a rational representation $V$ of $G$, 
we can form the associated vector bundle $E(V)$ on $X$.  By section 4
above, $E(V)$ acquires 
a flat connection $\nabla_V$.  Further, given two $G$ modules $V$ and $W$, and a 
$G$-module homomorphism $\varphi: V\rightarrow W$,  we obtain a homomorphism of vector  
bundles 
$\stackrel{\sim}{\varphi} : E(V)
\rightarrow E(W)$.   We note that $\stackrel{\sim}{\varphi}$ is equivariant for the 
connections $\nabla_V$ and $\nabla_W$.  In other words, the 
following diagram commutes
$$
\xymatrix{
E(V)\ar[d]_{\widetilde{\varphi}} \ar[rr]^{\nabla_V} && E(V) \otimes 
\Omega^{1}_X \ar[d]^{\widetilde{\varphi}\otimes id_X }\\
E(W) \ar[rr]^{\nabla_W} && E(W) \otimes \Omega^{1}_X }
$$

We now consider a faithful representation of $G$ on a vector space $V$ so that 
$G \subset SL(V)$ (such a representation exists since $G$ is assumed to be semisimple).  
We consider the induced vector bundle $E(V)$ and the flat connection $\nabla_V$ on $E(V)$.  
We define the following category

\begin{definition}
$\C(E,V) = \{$ vector bundles  $S$ on  $X$ with a flat connection  such
 that there  exist $A \subset 
B \subset\oplus_i E(V)^{\otimes n_i} \otimes E(V)^{*^{\otimes m_i}}$, with  $A$ and $B$ 
preserved 
by the connection on $\oplus E(V)^{\otimes n_i} \otimes E(V)^{*\otimes m_i}$, and $S =B/A$ with the induced flat connection, 
where the direct sums are finite direct sums $\}$. 
\end{definition}

We now recall

\begin{Theorem}
Let $\C$ be a category with a distinguished object 
$\O$ 
equipped 
with an operation
$$
\stackrel{\wedge}{\otimes} :\C \times \C \rightarrow \C
$$
and let $T: \C\rightarrow k$-mod be a  functor satisfying the following eight conditions:\\
\begin{enumerate}
\item $\C$ is an abelian category with direct sums;
\item isomorphism classes of objects of $\C$ form a set;
\item $T$ is an additive, faithful and exact fuctor;
\item $\stackrel{\wedge}{\otimes}$ is $k$-linear in each variable, and $T \circ
\stackrel{\wedge}{\otimes} =\otimes  \circ (T \times T)$;
\item $\stackrel{\wedge}{\otimes}$ is associative preserving $T$;
\item  $\stackrel{\wedge}{\otimes}$ is commutative preserving $T$;
\item the object $\O$ of $\C$ is equipped with an isomorphism $\phi: k \rightarrow T(\O)$ 
such that $\O$ is an identity object of $\stackrel{\wedge}{\otimes}$ preserving $T$;
\item  for every object $L$ of $\C$ such that $T(L)$ is one-dimensional,
there is an
object $L^{-1}$ of $\C$ such that $L \stackrel{\wedge}{\otimes} L^{-1}$ 
is 
isomorphic to $\O$.
\end{enumerate}

Then there exists a unique affine algebraic group scheme $H$ defined over $k$ 
(not necessarily of finite type) such that the quadruple $(\C, \stackrel{\wedge}{\otimes}, T, \O)$ is identified
with the quadruple $(H- \mbox{mod}, \ \stackrel{\wedge}{\otimes}, T, \O)$. 
A quadruple 
$(\C, \stackrel{\wedge}{\otimes}, T, \O)$ as in the theorem satisfying
the conditions is called a Tannakian category.
\end{Theorem}

\begin{Proof}See section 4, Chapter 2, pages 152--153 in \cite{S} and pages
118--120 in \cite{N}.\hfill{QED}
\end{Proof}

We now further recall the following.  Let $H$ be an affine algebraic group scheme over $k$ 
and $S$ a scheme over $k$.  Let $E_H\rightarrow S$ be a principal $H$-bundle on $S$.  
Then for every representation $\rho: H \rightarrow GL(V)$ in $H$-mod, we can construct the 
associated vector bundle $E_\rho =E_H(V)$.  This defines a functor $\F_E$ from $H$-mod to the 
category Vect $(S)$ of vector bundles on $S$.  We have

\begin{Theorem} 
Let $\F: H$-mod$\rightarrow \mbox{Vect}(S)$ be a 
functor satisfying the following:

\begin{enumerate}
\item $\F$ is a $k$- additive exact functor;
\item $\F\circ \stackrel{\wedge}{\otimes} =\otimes \circ(\F\times \F)$;
\item  $\F$ preserves commutativity, in other words, if $c$ is the canonical 
isomorphism of 
$V \stackrel{\wedge}{\otimes} W$ with $W\stackrel{\wedge}{\otimes}
V$ in $H$-mod, then $\F(c)$ is the canonical isomorphism of the corresponding vector 
bundles;
\item  $\F$ preserves associativity, i.e., if $a$ is the canonical isomorphism of $U  
\stackrel{\wedge}{\otimes} (V \stackrel{\wedge}{\otimes} W)$ with $(U
\stackrel{\wedge}{\otimes} V)\stackrel{\wedge}{\otimes} W$ in $H$-mod, then $\F(a)$ is the 
canonical isomorphism of the corresponding vector bundles; 
\item  the vector bundle $\F(\O)$ is the trivial line bundle $\O_S$ on $S$;
\item for any $ V\in H$-mod of dimension $r$, the vector bundle $\F(V)$ 
is of rank~$r$.
\end{enumerate}

Then there exists a unique principal $H$-bundle $E\rightarrow 
S$ such that 
$\F$ is identified with $\F_E$.
\end{Theorem}

\begin{Proof}See \cite{N}, Lemma 2.3, Proposition 2.5.\hfill{QED}
\end{Proof}

 We now observe that the category $\C(E,V)$ defined above, with 
 $\stackrel{\wedge}{\otimes}$ as tensor product of vector bundles, $T$ 
 as the functor which assigns to a vector bundle $W \in \ \C(E, V)$ the 
fibre $W_x$ of $W$ at a fixed point
 $x \in X, \ \O$ as the trivial line bundle $\O_X$ on $X$, 
and homomorphisms preserving connections, is a 
Tannakian category as in Theorem 5.1 above.  Hence there is a affine 
group scheme $M$ over $k$
 associated  to $\C(E,V)$,  which we call the monodromy, by Theorem 5.1 
above.

 We have an identification $M$-mod $=\C(E,V)$ where $M$-mod denotes the 
 category of rational $M$-modules as before.  Since every object of 
$\C(E,V)$ is a vector
 bundle on $X$, we obtain a natural functor 
$$\F:M-\ {\rm mod} \  \rightarrow  \ 
{\rm Vect} \  (X),$$
 and it is easy to see that $\F$ satisfies all the conditions of 
Theorem 5.2 above.  We thus obtain a $M$-bundle,
$$
\pi_1: E_M\rightarrow X \ \mbox{on}\ X.
$$

We now observe that every $G$-module is in a natural way an $M$-module 
as follows: Given a $G$-module $W$,  it can be obtained as a 
subquotient of a finite direct sum $\oplus_i V^{\otimes n_i}\otimes 
V^{*\otimes m_i}$ where $V$ is the faithful $G$-module we have chosen at 
the outset.  Hence the connection on $E(W)$ is preserved by the 
connection on $\bigoplus_{i} E(V)^{\otimes n_i}
\otimes  E(V){^*{^{\otimes m_i}}}$. Thus $E(W)$ is in the category ${\cal 
C}(E, V)$ and is hence an $M$-module.  We thus obtain a natural functor
$$G - {\rm mod }  \to M - \ {\rm mod}$$
and hence a homomorphism of group schemes $M \to G$.
We now show that every $M$-module is a subquotient of a $G$-module,
as follows: an $M$-module is an object of ${\cal C}(E, V)$, i.e., is
a vector bundle $W_1$ on $X$ with a flat connection which is a
subquotient of $\bigoplus_{i} E(V)^{\otimes n_i}
\otimes  E(V){^*{^{\otimes m_i}}}$ preserved by the flat connection.
Since $\pi^*E(V)$ is trivial, with the pullback connection 
$\pi^*(\nabla)$ being the exterior differential (see Appendix 2, 
Proposition B.1),
it follows that
$\pi^*(W_1)$ is also trivial with the pullback of the connection
being the exterior differential $d$. We thus obtain a vector space subquotient
$W$ of $\bigoplus_{i} V^{\otimes n_i}
\otimes  V{^*{^{\otimes m_i}}}$ with an $M$ action on $W$ induced by
the $G$ action on $V$ through the already defined homomorphism
$M \to G$. Since $\bigoplus_{i} V^{\otimes n_i}
\otimes  V{^*{^{\otimes m_i}}}$  is a G module, we have obtained an
$M$-module $W$ as a subquotient of a $G$-module such that $W_1$ is
isomorphic to the associated bundle $E_M(W)$ with the induced
flat connection. This proves our claim.

Since every $M$ 
module is a subquotient of a $G$-module, by Proposition 2.21 (b), page 
139, in Deligne-Milne's paper(see \cite{DM}), $M$ is a closed subgroup scheme
of
$G$.  
If $\F_E$ is the functor $\F_E: G$ mod $\to Vect (X)$ which associates 
to a rational $G$-module $W$ the vector bundle $E(W)$ on $X$, we obtain 
a commutative diagram 
$$\xymatrix{G - {\rm mod}\ar[dr]_{\F_E} \ar[rr] 
&& M-{\rm mod} \ar[dl]^{\F_{E_M}} && \\ 
& Vect(X)} 
$$
and this shows that the principal bundle $E_M$ induces the $G$-bundle 
$E$ on $X$, by Theorem 5.2 above.  We have thus obtained a reduction 
of structure group $E_M \subset E$ of $E$ to $M$.

\begin{Theorem}
Let $k$ be the field of complex numbers.  Let the topological fundamental group of $X$ be trivial. Then $E_M$ is trivial.
\end{Theorem}

\begin{Proof}We first observe that since $X$ is simply
connected, the universal cover of $X$ is algebraic, and the universal
covering map is algebraic, viz, the identity map $X\to X$. If
$\pi:E\to X$ is the $G$-torsor we started with, from our hypothesis
that $G$ is simply connected, and the homotopy exact sequence for
the topological fundamental group, we obtain that $E$ is also simply
connected. We observe that the use of the homotopy exact sequence
for an algebraic torsor can be justified as follows: we show
that $E \to X$ is a Serre fibration (see Definition 6.2, Chapter 7
in \cite{B}) and by Theorem 6.7, Chapter 7 in \cite{B}, we use the homotopy
exact sequence. To see that $E \to X$ is a Serre fibration, we observe
that by Theorem 6.11, Chapter 7, in \cite{B}, the property is local on
the base. Now,
$E \to X$ is locally isotrivial (i.e., every point has a finite etale
neighbourhood over which the torsor is trivial, see Remark 1.1 above,
and \cite{Ra}, Lemme XIV 1.4), and a finite etale cover is a 
Serre fibration since it is a topological covering space, while a
trivial bundle is clearly a fibration. It follows that $E \to X$
is a Serre fibration and we can apply the homotopy exact sequence. 
So the universal cover of $E$ is also algebraic, viz, the
identity map $E\to E$. Let $E_{SL(V)}$ denote the principal frame
bundle associated to the vector bundle $E(V)$. Since 
$\pi^*(E(V))$ is trivial on $E$,$\pi^*(E_{SL(V)})$ is the trivial
$SL(V)$ bundle on $E$.

      Let $\widetilde X$ denote the universal cover of $X$ (since
$X$ is simply connected, $\widetilde X =X$). Since $E(V)$ carries a
flat connection $\nabla$,$E_{SL(V)}$ is induced analytically by
a map $\widetilde X \to E_{SL(V)}$ inducing the flat structure.
Similarly, if $\widetilde E$ denotes the universal cover of $E$
($\widetilde E =E$ since $E$ is simply connected), the flat
connection  $\pi^*(\nabla)$ on $\pi^*E_{SL(V)}$ is induced 
analytically by a map $\widetilde E \to \pi^*E_{SL(V)}$
inducing the flat structure on $\pi^*E_{SL(V)}$.  
Further, by construction (see Section 4 above and Appendix 2, Proposition B.1),
the connection $\nabla$ on $E(V)$ has the property that 
$\pi^*(\nabla)=d$, the exterior differential, on the trivial
bundle $\pi^*(E(V))$. The map $\widetilde E \to \pi^*E_{SL(V)}$
is defined by a linearly independent set of sections $s_i$ of 
$\pi^*(E(V))$ such that $\pi^*(\nabla)(s_i)=0$. Since
$\pi^*(\nabla)=d$,$ds_i=0$, and so the $s_i$ are constant.
Since constant sections are algebraic,
it follows that the morphism 
$\widetilde E \to \pi^*E_{SL(V)}$ given by the flat connection (viz,$d$) on
$\pi^*(E_{SL(V)})$ is algebraic. From the commutative diagram
$$
\xymatrix{
E \ar[r] \ar[d]_{\pi} & \pi^* E_{SL(V)} \ar[d]\\
X\ar[r] & E_{SL(V)} }
$$
it follows that the morphism $X\to E_{SL(V)}$ given by the flat
connection $\nabla$ is also algebraic (we have to show that an
algebraic function on $E_{SL(V)}$ pulls back to an algebraic
function on $X$; however an algebraic function on $E_{SL(V)}$
pulls back to an algebraic function on $\pi^*(E_{SL(V)})$ since
$\pi$ is algebraic, and an algebraic function on $\pi^*(E_{SL(V)})$
pulls back to an algebraic function on $E$ since
$E \to \pi^*E_{SL(V)}$ is algebraic; since a function pulled back
from $E_{SL(V)}$ under $\pi$ is $G$-invariant we obtain a
$G$-invariant algebraic function on $E$ which gives an 
algebraic function on $X$).
We observe that the connection on every object of ${\cal C}(E, V)$
is induced from the connection on $E(V)$, and hence by the
construction of $M$ and $E_M$,the morphism
$X\to E_{SL(V)}$ factors through $E_M$. We further observe that
if $\phi$ denotes the map $\widetilde X \to E_{SL(V)}$
defining the flat structure, and $\widetilde \pi$ denotes the
projection $E_{SL(V)} \to X$,
then $\widetilde \pi \circ \phi$ is the
universal covering map $\widetilde X \to X$ (this can be seen
by noting that $E_{SL(V)}$ is defined by a representation of
the topological fundamental group of $X$ since it carries a
flat connection).
Thus we obtain an
algebraic section $X \to E_M$. Since a principal bundle with a
section is trivial, it follows that $E_M$ is trivial (algebraically).         
\hfill{QED}
\end{Proof}

\begin{Remark}
We observe that $\nabla$ has regular
singularities by construction, since $\pi^*(\nabla)=d$, and $d$ has
regular singularities. An alternative proof of the above theorem
can be given using Deligne's theorem(see \cite{D} and \cite{K}).
\end{Remark}

\section{The Main Theorem}

Let 
$X$ be a smooth, affine variety defined over $k$, where $k$ is an
algebraically closed field of characteristic $0$ contained in $\CC$.
Let $G$ be a 
semisimple algebraic group over $k$.  In this section, we prove 

\begin{Theorem} 
Let $X_{\CC}$ have trivial singular 
cohomology with integer coefficients (i.e. $H^0(X_{\CC}, \ZZ) = \ZZ$ and 
$H^i(X_{\CC}, \ZZ) = 0$ for $ i \ge 1$), let $\Omega^1_X$ be algebraically
trivial, and assume that $X_{\CC}$ is 
topologically simply connected.  Then any $G$-bundle $\pi: E \to X$ 
on $X$ over $k$ is trivial. (See Remark 6.2 and Remark 6.4).
\end{Theorem}

\begin{Proof} We first reduce to the simply 
connected case.  So we assume (for the purpose of this step), that the 
theorem is true for every semisimple, simply connected group over $k$.  
We now consider a semisimple group $G$.  If $G$ is not simply connected, 
then let $\widetilde{G} \to G$ be its simply connected cover.  We 
observe that the kernel $F = Ker (\widetilde{G} \to G)$ is central in 
$\widetilde{G}$.  We have 
$$1 \to F \to \widetilde{G} \to G \to 1$$
and the induced etale cohomology sequence
$$H^1(X, \widetilde{G})\to H^1 (X, G) \to H^2 (X, F)$$
We observe that we may assume without loss of generality that $F$
is the group $\ZZ/n$. Since we have assumed that $H^i(X_{\CC}, \ZZ)$
vanishes for $i=2,3$, we obtain $H^2(X,F)=0$ (we observe that
etale cohomology and singular cohomology with cyclic 
coefficients are isomorphic by Theorem 3.12, Chapter 3 in \cite{M}).

\medskip

Thus $H^1(X, \widetilde{G}) \to H^1(X, G)$ is surjective. 
 The $G$-bundle $E$ on $X$
can be lifted to a $\widetilde{G}$ bundle $\widetilde{E}$ on $X$.  Since 
$\widetilde{G}$ is semisimple and simply connected, if the bundle 
$\widetilde{E}$ is shown to be trivial, then $E$ is trivial as well. 
\medskip

Thus we have reduced to the case when $G$ is semisimple and simply 
connected. 
\medskip

For the rest of the proof, we assume that $G$ is semisimple and simply 
connected. 
\medskip

As before, let $\pi : E \to X$ be a $G$-bundle on $X$. 
Then by Lemma 2.1, we can assume that $k$ is the field of complex
numbers $\CC$. By Theorem 5.3 and the paragraph preceding it,
it follows that $\pi : E \to X$ is the trivial $G$-torsor.
\hfill{QED}
\end{Proof}

\begin{Remark} 
We observe that we have only used that
$H^i(X_{\CC}, \ZZ) = 0$ for $ i= 1,2,$ and $3$, $\Omega^1_X$ is 
algebraically trivial, and that $X_{\CC}$ is
topologically simply connected.
\end{Remark}

\begin{Remark} 
There are contractible, affine, smooth surfaces
over $\CC$ which are not isomorphic to affine two space(see \cite{GM}).
\end{Remark}

\begin{Remark} 
We note that when $X_{\CC}$ is contractible
(as is the case when it is simply connected and the singular cohomology
with integer coefficients is trivial), the tangent bundle of $X_{\CC}$
is topologically trivial, and since $X_{\CC}$ is Stein, this implies
that the tangent bundle is holomorphically trivial. This is enough
for the proof in Appendix 1 of Proposition A.1 to carry through, so Proposition 3.1
is valid under the hypothesis that $X_{\CC}$ is contractible.
The main theorem then states that if $X_{\CC}$ is contractible,
then every algebraic $G$-torsor on $X$ is algebraically trivial,
when $G$ is a semisimple algebraic group.
\end{Remark}

\appendix

\section{Appendix 1}

We have the splitting
$$
A^1_E =\Omega^1_X \oplus B^1_E
$$
with $B^1_E$ isomorphic to the Lie algebra bundle $E(\G)$.

There is natural inclusion of bundles
$$
\wedge^2 A^1 _E \ \subset A^2_E
$$
and thus we get the subbundle
$$
\Omega^2_X \oplus (\Omega^1_X \otimes B^1_E) \oplus \wedge^2 B^1_E \subset
A^2_E.
$$
The differential operator $d_E$ takes $H^0(X, B^1_E)$ to a $\CC$-vector
subspace $d_E H^0(X, B^1_E) \subset H^0(X,A^2_E)$

We consider the $\CC$- vector space spanned
$$
M =d_E H^0(X,B^1_E) +H^0(X,\Omega^1_X\otimes B^1_E)
$$
$M$ is a $\CC$-vector subspace of $H^0(X,A^2_E)$.

For $f\in H^0(X,\O_X)$ and $\eta \in H^0(X, B^1_E)$, we have,
$$
f \ d_E \eta =d_E (f\eta) -d\ f \ \otimes \ \eta
$$
Hence $M$ is a $H^0(X,\O_X)$ submodule of $H^0(X,A^2_E)$.

We now show 

\begin{proposition}
$M$ generates a vector subbundle of $A^2_E$ disjoint 
from $\Omega^2_X$. 
\end{proposition}

\begin{Proof}
Let $Y =E$ and $p:Y \to X$ be the projection $p = \pi$.
Then $p^* E =Y \times G$.

Let $\widetilde\pi   : p^* \ E \rightarrow Y$ so that
\[
\xymatrix{
p^*E \ar[d]_-{\widetilde{\pi}} \ar[r]^-{\widetilde{p}} & E \ar[d]^-\pi \\
Y \ar[r]_-p & X
}
\]


is a commutative diagram (fibre product).  Since $\pi^* \ B^1 _E \subset
\Omega^1 _E$, we have
$$
\widetilde\pi^* p^* B^1_E \subset  \Omega^1_{p^* E} =\Omega^1_Y \oplus
\Omega^1_G
$$
We also have 
$$
\pi^* d_E = d \ on \ \Omega^1_E
$$
so

$$
\widetilde\pi^*  d_E =\widetilde p^*  \pi^* \ d_E =
\widetilde\pi^*  d =d
$$
Further since $\pi^* \ B^1_E$ is isomorphic to $  E\times \G^*$ (the trivial
bundle),
we have
$\widetilde\pi^* p^* \ B^1_E$ is isomorphic to 
$Y \times G \times \G^*$ (the trivial bundle with fibre $\G^*$).\\
We have
$$
\begin{array}{ccc}
\widetilde\pi^*  p^* \ B^1_E  &\subset  \Omega^1_Y\oplus \Omega^1_G\\
&\downarrow_d& \\
&\Omega^2_Y \oplus (\Omega^1_Y\otimes \Omega^1_G) \oplus \Omega^2_G&\\
\end{array}
$$
Now, let $w_1,\ldots, w_n$ be $G$ invariant one forms on $G$ forming a global
frame (trivialisation on $G$) of $\Omega^1_G$.\\
(We observe that a $G$-invariant differential form on $G$ is nowhere vanishing
or identically zero).\\

We remark that $d w_i$ are $G$-invariant.\\
If $\sum_i \lambda_i d w_i (x)=0$ for some $ x \in G$ and $\lambda_i \in \CC$.\\
then
$$
\begin{array}{lll}
\sum \lambda_i d w_i(x) &=& d(\sum \lambda_i w_i) (x)\\
&=&0\\
\end{array}
$$
Since $\sum \lambda_i w_i$ is $G$-invariant, so is $d(\sum \lambda_i w_i)$ 
and since it vanishes at a point, it is identically zero.\\
Hence
$$
d(\sum \lambda_i w_i) =0 \ on \ G.
$$
Since $G$ is simply connected, $H^1(G, \CC) =0$ and by Grothendieck's algebraic
de Rham theorem,
$$
\sum \lambda_i w_i = df
$$
for a function $f$ on $G$.
Since $\sum \lambda_i w_i$ is $G$-invariant, for $ g \in G, \ g^* d\ f =df$
$$
\begin{array}{lllll}
\Rightarrow & d(g^* f) &=& df\\[1mm]
\Rightarrow & d(g^*  f- f) &=& 0\\[1mm]
\Rightarrow & g^* f -f &=& C_g
\end{array}
$$
where $C_g \in C$ is a constant depending on $ g \in G$.\\
We thus obtain a homomorphism $G\rightarrow G_a$ (where $G_a$ is the additive
group) by sending $g\mapsto C_g$.\\
Since $G$ is semisimple and $G_a$ is abelian, and there is no nonconstant
morphism from a semisimple algebraic group to an abelian group, we obtain
$C_g=0$ for all $ g \in G$.\\
Hence $f$ is $G$-invariant.  The only $G$-invariant function on $G$ is the
constant function, so $f$ is constant.  Therefore, $df=0$, hence
$$
\sum \lambda_i w_i =0
$$
contradicting the fact that the $w_i$ form  a frame at every point of $G$.
\\
Thus the $dw_i$ span a constant rank subbundle (in fact trivial) of 
$\Omega^2_{p^* E}$ disjoint from $\widetilde\pi^*  p^* \Omega^2_X$ and 
$\widetilde\pi^*  p^* (\Omega^1_X \otimes B^1_E)$.

The bundle $\widetilde\pi^*  p^* B^1_E$ is trivial.

Since $\Omega^1_X$ is trivial,
there are forms $\xi_1,\ldots, \xi_n \in H^0(X,\Omega^1_X)$ such that
$\xi_1 + w_1,\ldots, \xi_n +w_n$. form a global frame (trivialisation)
for $\widetilde\pi^* p^* B^1_E$. A $G$-invariant section of 
$\widetilde\pi^* p^* B^1_E$ is of the form 
$$
\sum_i f_i (\xi_i+w_i)
$$
where $ f_i \in H^0(Y, \O_Y)$. For a section $\eta$ of $B^1_E, \eta \in
 H^0(X, B^1_E), \widetilde\pi^*   p^* \eta$ is $G$-invariant so
$$
\widetilde\pi^*   p^* \ \eta =\sum_1 f_i (\xi_i + w_i)
$$
where $f_i \in H^0 (Y,\O_Y)$.\\
We now have
$$
\begin{array}{llll}
\widetilde\pi^*   p^*  d_E \eta & = d \widetilde\pi^*   p^* 
\eta\\[1mm]
&= d (\sum f_i (\xi_i +w_i))\\[1mm]
&= \sum_i df_i\otimes (\xi_i +w_i)\\[1mm]
& + \sum_i f_i(d \xi_i +dw_i)
\end{array}
$$
If at some point $ p =(y_0, g_0) \in Y \times G$ we have a linear relation
$$
\sum \lambda_i (d\xi_i (p) +dw_i(p))=0
$$
for $\lambda_i \in \CC$,\\
then since $ p =(y_0, g_0)$ and $d\xi_i$ is pulled back from $Y$  
while $dw_i$ is pulled back from $G$, we have
$$
\sum \lambda_i d \xi_i (y_0)=0
$$
and 
$$\sum \lambda_i \ d w_i(g_0)=0$$.\\
From the second relation above at $g_0$ and the earlier remark that the $dw_i$ are linearly independent at every point of $G$, it follows that
that $\lambda_i =0\ \forall i$.  Hence $d\xi_i+dw_i$ are linearly independent at
every point of $Y\times G$.\\
If
$$
\begin{array}{l}
\sum \lambda_i(df_i \otimes (\xi_i + w_i)) (P)
=\sum \mu_i f_i (d\xi_i +d w_i) (P)\\
\end{array}
$$
at some point $P$, for some $\lambda_i, \mu_i \in \CC$, then we have 
$$
\alpha_1 +\beta =\beta_2 +\gamma
$$
where
$$
\begin{array}{lllll}
\alpha_1 &=& \sum \lambda_i \ df_i \ \wedge\  \xi_i (P)\\[1mm]
&& \in q^*_1\ \Omega^2_{Y,P}\\[1mm]
\beta&=& \sum \lambda_i(df_i\otimes w_i) (P)\\[1mm]
&& \in q^*_1 \ \Omega^1_{Y,P} \otimes q^*_2 \ \Omega^1_{G,P}\\[1mm]
\alpha_2 &=&\sum  \mu_i \ f_i \ d \xi_i (P)\\[1mm]
&&\in  q^*_1 \ \Omega^2_{Y,P}\\
\end{array}
$$
and
$$
\begin{array}{lllll}
\gamma &=& \sum \mu_i f_i \ dw_i\\[1mm]
&& \in q^*_2 \ \Omega^2_{G,P}\\
\end{array}
$$
where $q_1, q_2$ are the two projections from $Y\times G$ to $Y$ and $G$
respectively.\\
From the direct sum decomposition 
$$
\begin{array}{llll}
q^*_1\Omega^2_Y &\oplus& (q^*_1\Omega^1_Y
\otimes q^*_2 \Omega^1_G) \\
&\oplus&  q^*_2 \Omega^2 _G
\end{array}
$$
it follows that $\beta(P)=0$ and $\gamma (P)=0$.\\
Now,
$$
\sum_I \mu_i f_i(P)  \ dw_i(P) =0
$$
implies that\\
$\mu_I=0$ whenever $f_i(P)\neq 0$.\\
Hence
$$
\sum \mu_I f_i(P) \ d\ \xi_i \ (P)=0
$$
since for every $i$, either 
$$
f_i(P) =0 \ or \ \mu_i=0
$$
Thus $\alpha_2(P)$ is also 0. So the intersection is 0.

We have thus shown that $M$ generates a vector subbundle of $A^2_E$ isomorphic
to
$$
(\Omega^1_X \otimes B^1_E)\oplus V
$$
where $V$ is a rank $n$ subbundle and this is disjoint from $\Omega^2_X$.
\end{Proof}

\section{Appendix 2}

On $E$, we have the following commutative diagrams:

$$
\begin{array}{ccccccc}
0 &&0 &\\
\downarrow && \downarrow\\
0\rightarrow \pi^* \Omega^1_X & \rightarrow &\pi^* A^1 _E &\rightarrow & 
Q_1 &\rightarrow &0\\
\downarrow && \downarrow&&\|\\
0\rightarrow \Omega^1_E &\rightarrow & A^1_{\pi^{*}E} &\rightarrow & Q_1 &
\rightarrow & 0\\ 
\downarrow && \downarrow\\
\Omega^1_{E/X} & = & \Omega^1_{E/X}\\
\downarrow && \downarrow\\

0 &&0 &\\
\end{array}
$$

and

$$
\begin{array}{ccccccc}
0 &&0 &\\
\downarrow && \downarrow\\
0 \rightarrow \pi^* \Omega_X^2 &\rightarrow & \pi^* A^2 _E &\rightarrow &
Q_2 &\rightarrow 0\\
\downarrow && \downarrow\ && \|\\
0 \rightarrow  \Omega_E^2 &\rightarrow &  A^2_{\pi^*E} &\rightarrow &
Q_2 &\rightarrow 0\\
\downarrow && \downarrow\\
\Omega^2_{E/X} & = & \Omega^2_{E/X}\\
\downarrow && \downarrow\\
0 &&0 &\\
\end{array}
$$
The splittings $s_1, s_2$ (see Proposition 3.1),
$$
\begin{array}{lll}
s_1 &:& A^1_E \rightarrow \Omega^1_X\\
s_2 &:& A^2_E \rightarrow \Omega^2_X\\
\end{array}
$$
chosen so that

$$ d_X \circ s_1 =s_2 \circ d_E
$$
are seen to induce splittings 
$$
\stackrel{\sim}{s_i} : A^i_{\pi^*E} \rightarrow \Omega^i_E \ i=1,2
$$
such that
$$
d \ \circ \ \stackrel{\sim}{s_1} =\stackrel{\sim}{s_2}\ \circ \ d_{\pi^*E}
$$
where $ d_{\pi^* E} : A^1_{\pi^*E} \rightarrow A^2 _{\pi^* E}$ is the 
operator defined for the bundle $\pi^* E$ on $E$ as before, and such that
$$
\stackrel{\sim}{s_i} \mid \pi^* A^i_E = \pi^* s_i \ i=1,2
$$
(This can be seen as follows: the splittings $s_i$ give
splittings $Q_i \to \pi^*A^i_E$ and these give splittings 
$Q_i \to A^i_{\pi^*E}$, and these give the $\stackrel{\sim}{s_i}$.)
Since $\pi^*E$ is canonically trivial (namely, by the isomorphism $E\times G \simeq E \times_X E$ by the map $(e,g) \mapsto (e,eg)$) on $E$, we have a canonical trivialisation $\pi^* E(V) = E 
\times V$ on $E$. Let $\pi^* (\nabla) = d+\omega$ on $\pi^*E(V) = E\times V$
where $\omega$ is a matrix of 1-forms on $E$.  We observe that $\pi^* (\nabla) 
(=d+\omega)$ on $ E \times V$ is $\stackrel{\sim}{s_1}\ \circ \ d_{\pi^* E}$ and
$\pi^*(\nabla) (=d+\omega)$ on $E \times V \otimes \Omega^1 _E$ is $\stackrel{\sim}
{s_2} \ \circ \ d_{\pi^* E} \ \circ \ i_1$, where $ i_1$ is the inclusion
$$
E\times V \otimes \Omega^1_E \subset  E \times V \otimes  A^1_{\pi^*E}
$$
We have the commutative diagram
$$
\begin{array}{cccccc}
E \times V &\stackrel{d_{\pi^* E}}{\longrightarrow} & E \times V \otimes 
A^1_{\pi^*E}  &\stackrel{d_{\pi^* E}}{\longrightarrow} &  E\times V \otimes A^2_{\pi^*E}\\
\| && \stackrel{\sim}{s_1} \downarrow \ \ \ \ \uparrow i_1 && 
\stackrel{\sim}{s_2}\downarrow \ \ \ \ \uparrow i_2\\
E\times V &\stackrel{\pi^*(\nabla)}{\longrightarrow} &  E\times V \otimes
\Omega^1_E&\stackrel{\pi^*(\nabla)=d+\omega}{\longrightarrow} &E\times 
V\otimes \Omega^2_E 
\end{array}
$$
Since $\pi^*(\nabla)=d+\omega$ from $E\times V\otimes \Omega^1_E$ to $E\times V 
\otimes \Omega^2_E$ is equal to $\stackrel{\sim}{s_2}\ \circ \ d_{\pi^*E} \ 
\circ \ i_1$, and by construction $d \ \circ \ \stackrel{\sim}{s_1} =
\stackrel{\sim}{s_2} \ \circ \ d_{\pi^*E}$, it follows that $\pi^*(\nabla) 
=d$ on $E\times V \otimes \Omega^1_E$.  In other words, $\omega=0$ and we have shown 

\begin{proposition}
$\pi^*(\nabla)=d$.
\end{proposition}


\noindent
School of Mathematics \\
Tata Institute of Fundamental Research \\
Colaba, Mumbai 400 005 \\
India\\
{\bf e-mail:subramnn@math.tifr.res.in}

\end{document}